\documentclass{amsart}
\usepackage{color}
\usepackage{graphicx}
\usepackage{amsmath,amssymb}
\theoremstyle{plain}

\newtheorem{thm}{Theorem}
\newtheorem*{thm*}{Theorem}

\newtheorem*{lem*}{Lemma}
\newtheorem*{cor*}{Corollary}
\newtheorem{cor}{Corollary}
\newtheorem*{conj*}{Conjecture}
\newtheorem*{prop*}{Proposition}
\newtheorem{prop}{Proposition}
\theoremstyle{definition}

\newtheorem{exam}{Example}
\newtheorem{defn}{Definition}
\newtheorem*{defn*}{Definition}

\DeclareMathOperator{\ecken}{vert}

\DeclareMathOperator{\sym}{Sym}

\title[\rm Borsuk, Ulam, Tucker and constructive combinatorics]{The
  Borsuk-Ulam-property, Tucker-property\\ and constructive proofs in
  combinatorics} \author{Mark de Longueville and Rade T.\ \v Zivaljevi\'{c}} \address{Mark de
  Longueville, Freie Universit\"at Berlin, Fachbereich Mathematik, WE
  II, Arnim\-allee 3--5, 14195 Berlin, Germany.}
  \email{delong@math.fu-berlin.de}
  \address{Rade T.\ \v Zivaljevi\'{c}, Mathematical Institute SANU, Knez Mihailova 35/1, p.f.\ 367,
  11001 Belgrade, Serbia and Montenegro.}
\email{rade@turing.mi.sanu.ac.yu}
\begin{document}

\thispagestyle{empty}
\begin{abstract}
  This article is concerned with a general scheme on how to obtain
  constructive proofs for combinatorial theorems that have topological
  proofs so far. To this end the combinatorial concept of
  Tucker-property of a finite group $G$ is introduced and its relation
  to the topological Borsuk-Ulam-property is discussed.  Applications
  of the Tucker-property in combinatorics are demonstrated.
\end{abstract}
\maketitle
\section{Introduction}
Topological combinatorialists prove combinatorial theorems by means of
topological tools. For many combinatorialists there seems to remain an
unsatisfactory feeling arising from the expectation that a purely
combinatorial question also deserves a combinatorial proof. In many
cases this is substantiated by the quest to find algorithms that
provide solutions \cite{alon-constructive,simmons-su,su-housemates}.

An outstanding application of topological methods and the proof
that offered a completely new perspective on topological
combinatorics was Lov\'asz' proof \cite{lovasz} of the Kneser
conjecture in 1978. It has an interesting history of
simplifications and generalizations, which culminated in a
combinatorial proof by Matou\v{s}ek \cite{matousek-kneserproof} in
2000. Other questions to be mentioned are various fair division
problems \cite{su-housemates,simmons-su}, one closely related to
the necklace problem \cite{alon-necklace}, and further coloring
problems of graphs and hypergraphs
\cite{matousek-ziegler,ziegler}.

At the heart of Matou\v{s}ek's combinatorial proof of Lov\'asz'
theorem is the deeper understanding of the relation of the
Borsuk-Ulam-theorem to its combinatorial counterpart, a lemma by
Tucker \cite{tucker,freund-todd}.

In this article we will investigate and generalize such a
correspondence and show its applicability towards the necklace
problem and the related fair division problem, or more generally
towards any problem solvable by the relatives of the
Borsuk-Ulam-theorem. This might eventually lead to a general
recipe to produce constructive, combinatorial proofs for theorems
in topological combinatorics.

\medskip {\bf Outline of the paper.} In order to introduce the setup
and to motivate the subsequent generalizations, we start with the
correspondence of the Borsuk-Ulam-theorem and the lemma by Tucker,
followed by a short discussion of Matou\v{s}ek's proof of Lov\'asz'
theorem. This provides a motivation for the introduction of
Borsuk-Ulam pairs and their combinatorial counterpart: Tucker triples.
In turn this leads to the Borsuk-Ulam-property for a group $G$, a
concept originally introduced by Sarkaria \cite{sarkaria}, and the
corresponding combinatorial property referred to as the
Tucker-property for $G$.  We investigate the connection between these
concepts and discuss its applications towards combinatorial problems.

\section{The Borsuk-Ulam-theorem and the Tucker-lemma}
\label{sec:review}

Let us recall the Borsuk-Ulam-theorem, a topological result which
is often illustrated to a layman by the claim that at any moment
there is a pair of antipodal points on the surface of the earth
with the same temperature and air pressure. This theorem has many
nice proofs usually applying some argument from homology theory
\cite{bredon}.
\begin{thm*}[Borsuk--Ulam]
  Let $f:\mathbb S^n\rightarrow \mathbb R^n$ be a continuous,
  antipodal map, i.e., a map such that $f(-x)=-f(x)$ for each $x\in
  \mathbb S^n$.  Then there exists an $x$ with $f(x)=0$.
\end{thm*}
As it turns out the Borsuk-Ulam-theorem can easily be derived from the
following combinatorial lemma, and vice versa.  An antipodally
symmetric triangulation of the sphere is a triangulation $K$ with the
property that $\sigma \in K$ implies $-\sigma \in K$.
\begin{lem*}[Tucker \cite{tucker}]
  Let $K$ be an antipodally symmetric triangulation of the
  $n$-sphere\, $\mathbb S^n$ refining the triangulation of the sphere
  induced by the coordinate hyperplanes. Let
  $\lambda:\ecken(K)\rightarrow \{\pm1,\ldots,\pm n\}$ be an antipodal
  labeling of the vertices of $K$, i.e., $\lambda(-x)=-\lambda(x)$ for
  all $x\in \ecken(K)$. Then there exists $i\in\{1,\ldots, n\}$, and
  an edge in $K$ whose vertices are labeled by complementary labels
  $-i$ and $+i$.
\end{lem*}
An elementary, constructive proof of the lemma was given by Freund and
Todd~\cite{freund-todd}. For the proof and its relation to the
Borsuk-Ulam-theorem we refer to Matou\v{s}ek's book
\cite{matousek-book}. The requirement that the triangulation ``refines
the triangulation of the sphere induced by the coordinate
hyperplanes'' can be weakened \cite{prescott-su}. It actually can be
removed, but then the proof proceeds by a detour using the continuous
Borsuk-Ulam-theorem, for details see \cite{matousek-book}. But in this
article we do not want to consider such detours.

We want to comment on the word ``constructive'': The proof of Tucker's
lemma by Freund and Todd is based on the construction of a particular
graph of degree at most two.  Following a path in this graph starting
at a known vertex of degree one, the inclined mathematician will end
up in a vertex corresponding to the desired edge. In order to do this
one actually will only need to construct the graph along this path. In
general this will be much quicker then to search all edges of the
triangulation. The reader might be reminded of the proof of the
Sperner lemma which indeed is very similar, and we point out that
Sperner's lemma is the combinatorial counterpart to Brouwer's fixed
point theorem \cite{aigner-ziegler}, a connection that has proved to
be very fruitful in finding combinatorial proofs and algorithms as
well (see e.g., \cite{su-housemates}).
\section{Matou\v{s}ek's proof of the Kneser conjecture and a general scheme}
\label{sec:matousek-proof}
For self-containment, we recall Kneser's original conjecture
\cite{kneser} from 1955 first proved by Lov\'asz in 1978: For every
partition of $\binom{[n]}{k}$ into $n-2k+1$ sets
$C_1,\ldots,C_{n-2k+1}$ there exists a set $C_i$ containing a pair of
disjoint $k$-sets.  
Here and in the sequel $[n]$ denotes the set
$\{1,\ldots,n\}$. $\binom{[n]}{k}$ denotes the $k$-subsets of $[n]$.

Ji\v{r}\'i Matou\v{s}ek found an ingenious combinatorial proof
\cite{matousek-kneserproof} of the Kneser conjecture in 2000. It is
based on the insight that the combinatorial counterpart to the
Borsuk-Ulam-theorem, namely Tucker's lemma, is only needed in a very
mild form.  He presents a two-step procedure in order to obtain the
direct proof.  First he shows how a special case of Tucker's lemma
suffices, i.e., the octahedral Tucker lemma. And in a second step he
eliminates all intermediate steps and provides a compact elegant
constructive proof of Lov\'asz' theorem. Although Matou\v{s}ek
presents this as a proof by contradiction, his proof in fact finds a
pair of disjoint $k$-sets constructively.

Matou\v{s}ek's proof suggests a general scheme to obtain constructive
proofs for theorems with topological proofs. Note that step
\textit{iii)} might be easier to achieve if it is known how the
combinatorial result implies the topological one.

\begin{enumerate}
\item[\textit{i)}] Find a combinatorial counterpart to the topological
  theorem used in the proof.
\item[\textit{ii)}] Identify a special case of the combinatorial
  statement as needed, and prove it directly.
\item[\textit{iii)}] Replace the topological argument and the
  accompanying spaces etc.\ by the combinatorial counterparts.
\end{enumerate}

\section{Borsuk-Ulam pairs and Tucker triples}
\label{sec:pairs-triples}

The relationship between the Borsuk-Ulam theorem and Tucker lemma,
reviewed in Section~\ref{sec:review}, deserves a closer analysis
and motivates the introduction of more general concepts.
For all topological and combinatorial concepts that appear we refer to
\cite{matousek-book}, see also \cite{elisat}.

Let $K$ be a finite simplicial complex with a simplicial action of
a finite group $G$. Let $V$ be a real representation of $G$ and
$Q\subset V$ a $G$-invariant convex polytope such that $0\in {\rm
int}(Q)$.

\begin{exam}
  As a special case of such a complex $K$ the following is playing a
  special role.  Let $G$ be a non-trivial finite group of order $k$,
  let $n\geq 1$, and let $N:=n(k-1)$. Consider $E_NG=G\ast\ldots\ast
  G$, the $(N+1)$-fold join of $G$ with itself, where $G$ is regarded
  as a $0$-dimensional simplicial complex. We donote the vertex set of
  this complex by $G\times [N+1]$. We will denote the elements of the
  (geometric realization of) $E_NG$ by $(t_0\cdot g_0,\ldots,t_N\cdot
  g_N)$ where $t_i\geq 0$, $\sum t_i=1$, $g_i\in G$.  $E_NG$ is a
  compact $N$-dimensional, $(N-1)$-connected space with a free
  $G$-action, given by the diagonal action $g\cdot(t_0\cdot
  g_0,\ldots,t_N\cdot g_N)=(t_0\cdot gg_0,\ldots,t_N\cdot gg_N)$.

\end{exam}

\begin{defn}
\label{def:BU-pair} A pair $(K,V)$ is called a {\em Borsuk-Ulam
  pair} for the group $G$, or just a Borsuk-Ulam pair if $G$ is fixed
in advance, if each $G$-equivariant map $f: K\rightarrow V$ has a
zero, that is if \,$0\in {\rm Image}(f)$.
\end{defn}

 \begin{defn}\label{def:Tucker-triple}
   A triple $(K,V,Q)$ is called a {\em Tucker triple} for a group $G$,
   or just a Tucker triple for short, if for each $G$-equivariant map
   (labelling) $\phi : \ecken(K)\rightarrow \ecken(Q)$, there exists a
   simplex $\sigma \in K$ such that $0\in {\rm
     conv}(\phi(\ecken(\sigma)))$.
 \end{defn}

 The Tucker lemma can be rephrased as the statement that
 $(K,\mathbb{R}^n, \Diamond^n)$ is a Tucker triple for the group
 $\mathbb{Z}_2$ where $\Diamond^n:={\rm conv}\{+e_i,-e_i\}_{i=1}^n$ is
 the crosspolytope in $\mathbb{R}^n$, and $K$ is a
 $\mathbb{Z}_2$-complex homeomorphic to $\mathbb{S}^{n}$ with the
 symmetric triangulation. The following example shows that such a
 ``crosspolytope'' exists for each finite group $G$.

\begin{exam}\label{exam:crosspolytope}
  {\rm Suppose that $G$ is a group of order $k$. Let $\mathbb{R}^k
    \cong {\rm span}\{e_g\mid g\in G\}\cong \mathbb{R}[G]$ be the real
    regular representation {\rm \cite{FulHar}} of $G$ and $W_G
    :=\{x\in \mathbb{R}^k\mid x_1+\ldots+ x_k=0\}$ the representation
    obtained by taking orthogonal complement of the diagonal, i.e.,
    the trivial $1$-dimensional representation. The generalized
    crosspolytope $\Diamond_k^n = \Diamond_G^n$ is defined as the
    convex hull of the union $\cup_{i=1}^n\Delta_{(i)}\subset (W_G)^n$
    where $\Delta_{(i)}$ is the simplex in the $i$-th copy
    $(W_G)_{(i)}\cong W_G$, spanned by the projections (of $\mathbb
    R^k$ onto $(W_G)_{(i)}$) of the orthonormal basis vectors
    $e_1,\ldots,e_k$. The polytope $\Diamond_G^n$ is clearly a
    $G$-invariant subspace of $(W_G)^n$ such that $0\in {\rm
      int}(\Diamond_G^n)$. Note that $\Diamond_G^n$ depends only on
    the order $k$ of the group $G$ which justifies the notation
    $\Diamond_k^n = \Diamond_G^n$.}
\end{exam}

Tucker triples are easily generated from the Borsuk-Ulam pairs. A
moment of reflection shows that each Borsuk-Ulam pair $(K,V)$ can be
upgraded to a Tucker triple $(K,V,Q)$ where $Q$ is an arbitrary
$G$-invariant polytope in $V$ such that $0\in {\rm int}(Q)$.  Indeed,
a $G$-equivariant labelling $\phi : \ecken(K)\rightarrow \ecken(Q)$ is
linearly (simplicially) extended to a $G$-equivariant map $f :
K\rightarrow V$ and a zero of $f$ is inside a simplex $\sigma$ such
that $0\in {\rm conv}(\phi(\ecken(\sigma)))$.

The converse is not true. The following proposition shows that
there may be striking differences between the two notions.

\begin{prop}\label{prop:Cara} 
  Assume that $G$ is a group of order $k$ and let $V$ be a real
  $G$-representation of dimension $N=n(k-1)$ such that $V^G=\{0\}$,
  i.e., such that $g\cdot x=0$ for all $g\in G$ if and only if $x=0$. Let
  $Q\subset V$ be a simplicial, $G$-invariant polytope such that $0\in
  {\rm int}(Q)$.  Then $(E_NG,V,Q)$ is always a Tucker triple for the
  group $G$.
\end{prop}

An example for a space $V$ as in the proposition is the representation
$(W_G)^n$ where $W_G$ is the $(k-1)$-dimensional, real
$G$-representation obtained from the regular representation by
factoring out the $1$-dimensional trivial representation.

\begin{proof}
  The result is an easy consequence of
the following remarkable result from convex geometry due to
I.~B\'{a}r\'{a}ny \cite{Barany, Barany-Onn}.

Suppose that
\begin{equation}\label{eqn:matrix}
\Omega = \left[\begin{array}{cccc}
v_{1,1} & v_{1,2} & \dots & v_{1,d+1}\\
v_{2,1} & v_{2,2} & \dots & v_{2,d+1}\\
\vdots & \vdots & \ddots & \vdots\\
v_{m,1} & v_{m,2} & \dots & v_{m,d+1}
\end{array}\right]
\end{equation}
is a matrix where the entries $v_{i,j}$ are vectors in a vector
space $\mathbb{R}^d$. Moreover, we assume that $0\in{\rm
conv}\{v_{i,\nu}\}_{i=1}^m$ for each $\nu$, i.e., that the origin
is in the convex hull of each column of the matrix $\Omega$. Then
there exists a function $\alpha : [d+1]\rightarrow [m]$
such that
\[
0\in {\rm
conv}\{v_{\alpha(1),1},v_{\alpha(2),2},\ldots,v_{\alpha(d+1),d+1}\}.
\]
Suppose that $\phi : \ecken(E_NG)\rightarrow \ecken(Q)$ is a
$G$-equivariant labelling. Assume $m:=k$ and $d:=n(k-1)=N$.  Let
$\Omega = [v_{g,j}]$ be the vector-valued matrix defined by $v_{g,j}:=
\phi((g,j))$ for each $(g,j)\in \ecken(E_NG)= G\times[N+1]$. For
each $\nu$,
\[
x_\nu:=\sum_{g\in G}v_{g,\nu} = \sum_{g\in G} \phi(g(e,\nu)) =
\sum_{g\in G} g\phi(e,\nu)
\]
is a $G$-invariant element in $V$. By the assumption $V^G=\{0\}$,
hence $x_\nu=0$ for each $\nu$, and we conclude that the matrix
$\Omega$ satisfies the conditions of B\'{a}r\'{a}ny's theorem.
Consequently, there exists a function $\alpha : [N+1]\rightarrow
G$ such that $0\in {\rm conv}(\{v_{\alpha(i),i}\}_{i=1}^{N+1})$, or
equivalently,
\[
0\in {\rm conv}(\{\phi((\alpha(i),i))\}_{i=1}^{N+1}) = {\rm
  conv}(\phi(\ecken(\sigma)))
\]
where $\sigma\in E_NG$ is the simplex determined by the
function $\alpha$. 
\end{proof}

\bigskip There exist examples of groups $G$ of order $k$ such that
$(E_NG, (W_G)^n)$ is not a Borsuk-Ulam pair. Such an example is
provided already by the group $\mathbb{Z}_6$, see \cite{Ziv} part (II)
where it was shown that there exists a $\mathbb{Z}_6$-equivariant map
$f : \mathbb S^5\rightarrow \mathbb S^4\subset\mathbb{R}^5\cong
W_{\mathbb{Z}_6}$. Let $G=\mathbb Z_6$, i.e., $k=6$, and $n:=1$, hence
$N=5$. In this case the map $f$ also implies the existence of a
$\mathbb{Z}_6$-equivariant map $E_NG \rightarrow W_G\cong
\mathbb{R}^5$ without zeroes. So among the consequences of
Proposition~\ref{prop:Cara} is the following observation.

\begin{cor}
  There exists a finite group $G$ and a Tucker triple $(K,V,Q)$ such
  that $(K,V)$ is not a Borsuk-Ulam pair. For example one can take
  $G=\mathbb{Z}_6$, define $V:=W_{\mathbb{Z}_6}\cong \mathbb{R}^5$ as
  the real representation of $\mathbb{Z}_6$ obtained from the regular
  representation modulo the $1$-dimensional trivial representation,
  and choose $Q$ to be a $\mathbb{Z}_6$-invariant, $5$-dimensional
  simplex in $W_{\mathbb{Z}_6}$.
\end{cor}

\section{The Borsuk-Ulam-property of $G$}
\label{sec:bu-property} In order to pursue the scheme outlined in
Section~\ref{sec:matousek-proof}, we  first state a family of
generalizations of the Borsuk-Ulam-theorem that has proven to be
useful in topological combinatorics.

Consider the space $\mathbb{E}_{n,k}$ of all $(n\times k)$-matrices 
with real entries and the property that all row sums are zero. In
other words $\mathbb{E}_{n,k}$ is the space of all matrices orthogonal
to the space of all matrices with entries in each row being identical.
In particular, $\mathbb{E}_{n,k}$ has dimension $n(k-1)$.  By labeling
the columns with elements of $G$, $G$ acts on $\mathbb{E}_{n,k}$ by
column permutations, and the only fixed point of this action is the
zero matrix.

The reader can easily convince herself that $\mathbb{E}_{n,k} \cong
(W_G)^n$, where $W_G$ is the representation described in
Example~\ref{exam:crosspolytope}. In other words $\mathbb{E}_{n,k}$ is
just a more concrete presentation of the representation $W_G$. Note
that in particular, the simplex $\Delta_{(i)}\subset (W_G)_{(i)}$
corresponds in $\mathbb{E}_{n,k}$ up to scaling to the convex hull of
the set of matrices
\begin{align*}
  {\small
    \begin{bmatrix}
      0\\
      \vdots\\
      0\\
      (e_g-\frac{1}{k}\sum_{h\in G}e_h)^t\\
      0\\
      \vdots\\
      0
    \end{bmatrix}
  }
  \text{, }{g\in G}\text{,}
\end{align*}
where the non-zero entries are in row $i$

In the spirit of Sarkaria \cite{sarkaria} we introduce the
following definition.

\begin{defn}
A group $G$ of order $k$ has the \emph{Borsuk-Ulam-property} if
$(E_NG, \mathbb{E}_{n,k})$ is a Borsuk-Ulam pair for each $n\geq
1$. In other words for each $n\geq 1$, every $G$-equivariant
continuous map $f:E_NG\rightarrow \mathbb E_{n,k}$ must have a
zero.
\end{defn}

Let us briefly review the case of the group $G=\mathbb Z_2$. In
this case $E_NG=G^{*(n+1)}\cong\mathbb S^n$, and the $G$-action is
given by the antipodal map. The space $\mathbb{E}_{n,2}$ can be
identified with $\mathbb R^n$ together with the action $x\mapsto
-x$. We conclude that $G=\mathbb Z_2$ has the Borsuk-Ulam-property
which is just a restatement of the Borsuk-Ulam-theorem.

The following theorem is very important tool in topological
combinatorics with numerous and diverse applications.

\begin{thm}[\"Ozaydin \cite{ozaydin}, Sarkaria \cite{sarkaria}, 
  Volovikov \cite{volovikov}]\label{thm:prime-power} For $p$ prime,
  $r\geq 1$, the group $G=(\mathbb Z_p)^r$ has the
  Borsuk-Ulam-property.
\end{thm}

The groups $G=(\mathbb Z_p)^r$ are the only groups for which we know
that the Borsuk-Ulam-property holds. For more information about this
and related problems we refer the reader to \cite{bartsch}.\par

The previous theorem has been used in the proofs of numerous
combinatorial theorems, most notably the topological Tverberg
theorem and its relatives
\cite{matousek-book,sarkaria,Ziv,elisat}. In order to demonstrate
its strength and as an overture to the proof of
Theorem~\ref{thm:approx}, we present a short proof of a theorem by
Alon on simultaneous equipartitions (splitting) of a set of
probability measures (necklaces). Compared to Alon's original
approach and other existing  proofs, see  \cite{matousek-book} and
\cite{elisat} for the references, the proof doesn't offer new
ideas. However the exposition is smooth and short providing an
excellent example of how a zero of a continuous map encodes all the
information needed for the solution (equipartition) of a geometric
problem.

\begin{thm*}[Alon \cite{alon-necklace}]
  Let $\mu_1,\ldots,\mu_n$ be continuous probability measures on the
  unit interval and $k\geq 2$. Then it is possible to cut the interval
  in $n(k-1)$ places and to partition the $n(k-1)+1$ resulting
  intervals into $k$ families $F_1,\ldots F_k$ such that $\mu_i(\cup
  F_j)=\frac{1}{k}$ for all $i$ and $j$.
\end{thm*}
It is easy to see that the number of cuts is best possible in general.
\begin{proof}
  It is a straightforward combinatorial exercise to reduce the problem
  first to the case $k=p$ a prime number.  Then the elements of
  $E_NG$, with $G=\mathbb Z_p$ and $N:=n(p-1)$, define $n(p-1)$ cuts
  of the unit interval together with a partition $F_1,\ldots, F_p$ of
  the resulting intervals: consider $(t_0\cdot g_0,\ldots, t_N\cdot
  g_N)\in E_NG$ and define $x_{-1}:=0$ and $x_j:=\sum_{i=0}^j t_i$ for
  $j=0,\ldots, N$. Then the cuts are given by $x_0,\ldots,x_{N-1}$ and
  the resulting intervals are partitioned by setting
  \mbox{$F_i:=\{[x_{j-1},x_j]:j\in\{0,\ldots,N\}, g_j=i+p\mathbb
    Z\}$.} (Degenerate intervals with $x_{j-1}=x_j$ may be put into
  any of the $F_i$ since they have measure zero.) Next we define a map
  \begin{align*}
    E_NG&\longrightarrow \mathbb E_{n,p}\\
    (t_0\cdot g_0,\ldots, t_N\cdot g_N)&\longmapsto
    (E_{ij})_{i=1,\ldots,n\atop j=1,\ldots,p}
  \end{align*}
  where $E_{ij}:=\mu_i(\cup F_j)-\mu_i(\cup F_{j-1})$ with the
  $j$-indices considered modulo $p$. Note that $(E_{ij})$ has row sums
  equal to zero by construction.  With the columns of the matrix
  labeled appropriately by the elements of $G$ this map is continuous
  and $G$-equivariant. Hence by the previous theorem there exists a
  zero, which yields the desired cuts and the partition.
\end{proof}

\section{The Tucker-property for $G$}
\label{sec:tuck-prop}

The combinatorial counterpart to a Borsuk-Ulam pair is a Tucker
triple, as discussed in Section~\ref{sec:pairs-triples}. Similarly,
the Borsuk-Ulam-property for $G$ has a combinatorial counterpart,
referred to as the Tucker property for $G$.

Another motivation for introducing this concept comes from the
conjecture of Simmons and Su \cite{simmons-su}, discussed in
Section \ref{simmons-su-conj}. Since a Borsuk-Ulam pair can be
upgraded to a Tucker triple in many ways, depending on the choice
of a $G$-invariant polytope $Q$, it is clear that there does not
exist a unique way of defining the ``Tucker property'' for $G$.
For example a different generalization of Tucker's lemma has been
used by Ziegler in \cite{ziegler}. Our definition is based on the
choice of the generalized crosspolytope $\Diamond_k^n =
\Diamond_G^n$ introduced in Example~\ref{exam:crosspolytope}.

\medskip
As before, let $G$ be a non-trivial finite group of order $k$, and
$N=n(k-1)$ where $n\geq 1$. Consider a $G$-invariant triangulation
$K$ of $E_NG$, i.e., if $g\in G$ and $\sigma\in K$ then
$g\cdot\sigma\in K$. Furthermore, we assume that $K$ refines the
natural triangulation of $E_NG$ induced from the $(N+1)$-fold join
operation of the $0$-dimensional complex~$G$.

\begin{defn}\label{def:Tucker-property}
The group $G$ has the {\em Tucker property} if $(K,
(W_G)^n,\Diamond_k^n)$ is a Tucker triple for $G$ for each $n\geq
1$ and each $G$-invariant subdivision $K$ of the complex $E_NG$,
$N=n(k-1)$.
\end{defn}

A slightly more combinatorial reformulation of
Definition~\ref{def:Tucker-property} is the following. A group $G$ is said
to have the Tucker property if for all $n\geq 1$ and all $K$ as above,
every $G$-equivariant labelling $\lambda:\ecken(K)\rightarrow G\times
[n]$, i.e., a labelling such that $\lambda(g\cdot
v)=g\cdot\lambda(v)$, has the property that there exists an $i\in [n]$
and a $(k-1)$-simplex in $K$ whose vertices are labelled by
$\{(g,i):g\in G\}$. The equivalence of the two formulations follows
from the observation that if $0\in {\rm conv}(S)$ for some $S\subset
\ecken(\Diamond^n_k)$, then $\ecken(\Delta_{(i)})\subset S$ for some
$i$.

Let us consider the case of $G=\mathbb Z_2$ again. In this case,
$K$ turns out to be an antipodally symmetric triangulation of the
$n$-sphere refining the triangulation induced by the coordinate
hyperplanes. Hence $G=\mathbb Z_2$ has the Tucker-property by
Tucker's lemma.

More generally, a consequence of Theorem~\ref{thm:prime-power} and
Proposition~\ref{prop:BU2Tuck} is that the group $G=(\mathbb
Z_p)^r$ has the Tucker-property for each prime $p$ and arbitrary
$r\geq 1$.

\section{Borsuk-Ulam vs.\ Tucker-property}
\label{sec:bu-tuck-prop}

As already observed in Section~\ref{sec:pairs-triples}, if $(K,V)$
is a Borsuk-Ulam pair then $(K,V,Q)$ is a Tucker triple for any
$G$-invariant convex polytope in $V$. The following proposition is
an easy consequence.

\begin{prop}\label{prop:BU2Tuck}
If a group $G$ has the Borsuk-Ulam-property than it also has the
Tucker-property.
\end{prop}
\begin{proof}
Let $k=|G|\geq 2$, $n\geq 1$, $N=n(k-1)$ and let $K$ be a
$G$-invariant triangulation of $E_NG$ refining the natural
triangulation. Furthermore, let 
%
%
$\lambda : \ecken(K)\rightarrow \ecken(\Diamond^n_k)$ be a
$G$-equivariant map. Let $\Lambda : K \rightarrow \Diamond^n_k\subset
\mathbb E_{n,k}$ be the linear (affine) extension of this map. Since
by assumption $G$ has the Borsuk-Ulam property, there exists a simplex
$\sigma\in K$ and $x\in\sigma$, such that $\Lambda(x)=0$. It follows
that $0\in{\rm conv}(\Lambda(\ecken(\sigma)))$.
\end{proof}

The following theorem shows that the converse to
Proposition~\ref{prop:BU2Tuck} is also true.

\begin{thm}\label{thm:Tuck-BU}
If a group $G$ has the Tucker property than it also has the
Borsuk-Ulam property.
\end{thm}

\begin{proof} For the sake of contradiction assume that $(E_NG,
(W_G)^n)$ is not a Borsuk-Ulam pair for some $n\geq 1$. In other
words we assume that there exists a $G$-equivariant map $f :
E_NG\rightarrow  (W_G)^n$ such that $0\notin {\rm
Image}(f)$. By compactness we can assume that ${\rm
Image}(f)\subset \mathbb{R}^k\setminus U$ for some neighborhood of
$0$ and by rescaling we can assume that $U=\Diamond_k^n$. The
radial projection $R : \mathbb{R}^k\setminus \Diamond_k^n
\rightarrow \partial(\Diamond_k^n)$ to the boundary of the
crosspolytope is $G$-equivariant, so we can assume that ${\rm
Image}(f)\subset\partial(\Diamond_k^n)$. By the (equivariant) {\em
simplicial-approximation theorem}, there is a $G$-invariant
subdivision $K$ of the complex $E_NG$ and a simplicial,
$G$-equivariant map $g : K\rightarrow
\partial(\Diamond_k^n)$ approximating $f$ in a suitable sense.
Then the restriction of $g$ on the $0$-skeleton
$K^{(0)}=\ecken(K)$ defines a $G$-equivariant labelling $\phi :
\ecken(K)\rightarrow \ecken(\Diamond_k^n)$ which contradicts the
assumption that $(K, (W_G)^n,\Diamond_k^n)$ is a Tucker triple.
Indeed, $0\notin {\rm conv}(\phi(\ecken(\sigma)))$ for each
$\sigma\in K$ is a consequence of the fact that ${\rm
conv}(\phi(\ecken(\sigma)))\subset
g(\sigma)\subset\partial(\Diamond_k^n)$.
\end{proof}

\section{$G$\,-Tucker-property and combinatorial proofs}
\label{simmons-su-conj} In this section we  discuss a possible
application of the $G$-Tucker-property towards a combinatorial
problem: finding approximate solutions for the
concensus--$\frac{1}{k}$-division problem. This also relates to a
conjecture of Simmons and Su, which we will discuss as well.\par
\noindent {\bf The necklace problem and
concensus-$\frac{1}{k}$-division.}  In \cite{alon-necklace}, Alon
investigates the theft of a necklace with beads of $n$ different
colors by $k$ thieves.  Under the assumption that there are a
multiple of $k$ beads of each color and that the necklace is
opened at the clasp, the question is whether it is possible to
always cut the necklace at $n(k-1)$ places and to distribute the
resulting pieces among the thieves in such a way that each of them
gets the same number of beads of each color.  In order to show
that this is indeed the fact, Alon proved the equipartition of
measures theorem that we discussed in Section
\ref{sec:bu-property} and showed how it applies to the discrete
situation.\par

In \cite{simmons-su}, Simmons and Su consider the question of
subdividing an object into two portions in such a way that $n$ given
people believe that the two portions are equal in value. If the
problem is modeled in terms of simultaneously equipartitioning a set
of $n$ measures, the existence of such a partition is given by Alon's
theorem for $k=2$, but there is no algorithm on how to obtain such a
solution. Using Tucker's lemma Simmons and Su desribe an algorithm to
obtain an $\varepsilon$-approximate solution to the problem for any
given $\varepsilon>0$. Here we want to address the generalization of
this problem already mentioned in \cite{simmons-su}: Subdividing an
object into $k$ portions such that according to the $n$ individual
measures all the portions have value $\frac{1}{k}$.  Again the
existence of such a subdivision is guaranteed by Alon's theorem.  But
what about algorithmic $\varepsilon$-approximations? \par

This problem can be divided into two steps: Finding a constructive
proof of the Tucker-property for $G$, and showing how the
Tucker-property can be applied to yield approximate solutions.  So far
we were only able to provide the second step, which we will
demonstrate here. As in the proof of Alon's theorem the approximation
problem easily reduces to the case $k=p$, $p$ prime.

\begin{thm}\label{thm:approx}
  Let $k=p$ be a prime, $\mu_1,\ldots,\mu_n$ be continuous probability
  measures on the unit interval, and $\varepsilon>0$ be given. Then a
  single application of the Tucker-property for $G=\mathbb Z_p$ yields
  $n(p-1)$ cuts of the unit interval together with a partition
  $F_1,\ldots,F_p$ of the resulting intervals, such that $|\mu_i(\cup
  F_j)-\frac{1}{k}|<\varepsilon$ for all $i$ and $j$.
\end{thm}
The following proof relies on ideas from \cite{simmons-su}, but has to
deal with some technical problems that do not occur in the case $k=2$.

\begin{proof}
  As in the proof of Alon's theorem in Section \ref{sec:bu-property},
  the elements $v$ of $E_NG$, $N:=n(p-1)$, encode $n(p-1)$ cuts of the
  unit interval together with a partition $\mathcal
  F(v)=\{F_1,\ldots,F_p\}$ of the resulting intervals.  By continuity
  of the $\mu_i$, let $K$ be a $G$-invariant triangulation of $E_NG$
  refining the natural triangulation with the property that for all
  pairs of neighboring vertices $v,w$ of $K$ with corresponding
  partitions $\mathcal F(v)=\{F_1,\ldots,F_p\}$ and $\mathcal
  F(w)=\{F_1',\ldots,F_p'\}$ the inequality
  \begin{align*}
    |\mu_i(\cup F_j)-\mu_i(\cup F_j')|<\frac{\varepsilon}{(p-1)^2}
  \end{align*}
  holds for all $i$ and $j$.  Again, such a triangulation can be
  obtained by iterated barycentric subdivision. We will now define a
  labeling
  \begin{align*}
    \lambda:\ecken(K)&\longrightarrow  G\times [n]\\
    v&\longmapsto(\lambda_1(v),\lambda_2(v)),
  \end{align*}
  where $\lambda_1$ and $\lambda_2$ are defined as follows. Let
  $v\in\ecken(K)$ with corresponding partition $\mathcal
  F(v)=\{F_1,\ldots,F_p\}$. Consider $m(v):=\min_{i,j} \{\mu_i(\cup
  F_j)\}$, and let $\lambda_2(v)$ be the smallest $i$ such that there
  exists a $j$ with $\mu_i(\cup F_j)=m(v)$. In order to define
  $\lambda_1(v)$ consider the sign vector
  $(\varepsilon_1,\ldots,\varepsilon_p)\in\{+,-,0\}^p$ defined by
 \begin{align*}
   \varepsilon_j:=
   \begin{cases}
     +,\quad \text{ if }\mu_{\lambda_2(v)}(\cup F_{j+1})>\mu_{\lambda_2(v)}(\cup F_j),\\
     -,\quad \text{ if }\mu_{\lambda_2(v)}(\cup F_{j+1})<\mu_{\lambda_2(v)}(\cup F_j),\\
     0,\quad\; \text{ if }\mu_{\lambda_2(v)}(\cup F_{j+1})=\mu_{\lambda_2(v)}(\cup F_j),
   \end{cases}
 \end{align*}
 with $j$-indices considered modulo $p$. If
 $(\varepsilon_1,\ldots,\varepsilon_p)=(0,\ldots,0)$ then
 $\mu_{\lambda_2(v)}(\cup F_1)=\mu_{\lambda_2(v)}(\cup
 F_2)=\cdots=\mu_{\lambda_2(v)}(\cup F_p)=\frac{1}{p}$ and hence, by
 definition of $\lambda_2(v)$, we have $\mu_i(\cup F_j)=\frac{1}{p}$
 for all $i$ and $j$. In this lucky event we found what we were
 looking for. \par

 If $(\varepsilon_1,\ldots,\varepsilon_p)\not=(0,\ldots,0)$, then in
 particular the sign vector is not constant, and we can define
 $\lambda_1(v):=[j]$, where $j\in\{1,\ldots,p\}$ is such that
 \begin{align*}
   (\varepsilon_j,\varepsilon_{j+1},\ldots,\varepsilon_p,
   \varepsilon_1, \ldots,\varepsilon_{j-1})
 \end{align*}
 is the lexicographic smallest vector among all cyclic permutations of
 the $\varepsilon$-vector, with respect to the linear order $-<0<+$ of
 $\{+,-,0\}$.  Thus we have defined a $G$-equivariant labeling. We can
 think of this labeling as saying which person $\lambda_2(v)$ is
 distressed most by the fact that according to its measure the portion
 $\cup F_{\lambda_1(v)}$ is the smallest with respect to all portions
 and measures.  By the Tucker-property for $G$, we obtain an $i_0\in
 [n]$ and a $(p-1)$-simplex $\{v_1,\ldots,v_p\}$ such that
 $\lambda(v_r)=([r],i_0)$ for $r=1,\ldots,p$. In other words, for
 person $i_0$ there exist $p$ different partitions very close to each
 other, such that the person is distressed about a different portion
 every time. But this means that they all must have similar size close
 to $\frac{1}{p}$.  Since person $i_0$ was most distressed, the
 portions must be similar in size for all other people as well. From
 here on we will just carry out the filthy details of this
 consideration.\par

 Let $\mathcal F(v_r)=:\{F_1^r,\ldots,F_p^r\}$, and define
 $x^{i}_{rj}:= \mu_i(\cup F^r_j)$.  We have the following properties.
 \begin{enumerate}
 \item For all $i,j$ and $r\not=r'$:
   $|x^{i}_{rj}-x^{i}_{r'j}|<\frac{\varepsilon}{(p-1)^2}$ by
   definition of the triangulation,
 \item for all $i,r$: $\sum_{j=1}^p x^{i}_{rj}=1$ since the $\mu_i$
   are probability measures,
 \item and we have for all $j$: $m(v_r)=x^{i_0}_{rr}\leq x^{i_0}_{rj}$
   since $\lambda(v_r)=([r],i_0)$.
 \end{enumerate}
 First, we will be concerned with the numbers $x^{i_0}_{rj}$.
 Properties (2) and (3) yield for all $r$: $x^{i_0}_{rr}\leq
 \frac{1}{p}$. Hence by (1), we obtain for all $r$ and $j$:
 $x^{i_0}_{rj}<\frac{1}{p}+\frac{\varepsilon}{(p-1)^2}$, which
 together with (2) yields
 $x^{i_0}_{rj}>\frac{1}{p}-\frac{\varepsilon}{(p-1)}$. Now let $i,j$,
 and $r$ be arbitrary. Then by the definition of the labeling
 $x^{i}_{rj}\geq
 m(v_r)=x^{i_0}_{rj}>\frac{1}{p}-\frac{\varepsilon}{(p-1)}$. By (2) we
 therefore obtain $x^{i}_{rj}<\frac{1}{p}+\varepsilon$. The last two
 inequalities yield the result.
 \end{proof}
 Note that the previous theorem also yields a proof of Alon's theorem
 by compactness of the space $E_NG$.\vspace{\baselineskip}

\noindent {\bf A conjecture of Simmons and Su.}
In \cite{simmons-su} Simmons and Su consider the following space
\begin{align*}
  S^n_k=\left\{(z_0,\ldots,z_n)\in\mathbb C^{n+1}:\sum_{i=0}^n|z_i|=1,z_i^k=|z_i|^k\right\}.
\end{align*}
Let $\omega:=e^\frac{2\pi i}{k}$, then the elements of $S^n_k$ have
the form $(t_0\omega^{j_0},\ldots,t_n\omega^{j_n})$ for some
$j_i\in\{1,\ldots,k\}$ and $t_i\geq 0$ with $\sum t_i=1$. The
symmetric group $\sym(k)$ acts on $S^n_k$ by
\begin{align*}
  \pi\cdot(t_0\omega^{j_0},\ldots,t_n\omega^{j_n})=
  (t_0\omega^{\pi(j_0)},\ldots,t_n\omega^{\pi(j_n)}).
\end{align*}
\begin{conj*}[Simmons \& Su \cite{simmons-su}]
  Suppose that $S^{n(k-1)}_k$ is triangulated invariantly with respect
  to the action of the symmetric group $\sym(k)$, and suppose that the
  vertices $V$ of the triangulation are labeled by a function
  \begin{align*}
    \ell:V\rightarrow \{\omega^jm:1\leq j\leq k,\; 1\leq m\leq n\},
  \end{align*}
such that for $\pi\in\sym(k)$ the condition $\ell(\pi(v))=\pi(\ell(v))$
holds for all $v\in V$. Then there must exist $k$ adjacent vertices in
the triangulation with labels \mbox{$\{\omega^jm:1\leq j\leq k\}$} for a
fixed $m$.
\end{conj*}
We will see that this conjecture holds in the case where $k$ is a
prime power $p^r$, in which case the requirement that everything is
symmetric with respect to the whole symmetric group can be weakened to
symmetry with respect to the subgroup $G=(\mathbb Z_p)^r$.
\begin{prop}
  Let $G$ be any group of order $k$. Then $G$ considered as a subgroup
  of $\sym(k)$ via an enumeration $\{g_1,\ldots,g_k\}$ of $G$ acts on
  $S^n_k$ and there is an equivariant homeomorphism from $S^n_k$ to
  $E_nG$.
\end{prop}
\begin{proof}
The homeomorphism is given by
$(t_0\omega^{j_0},\ldots,t_n\omega^{j_n})\longmapsto
  (t_0g_{j_0},\ldots,t_ng_{j_n})$.
\end{proof}
\begin{cor}
  The conjecture by Simmons and Su holds in the case $k=p^r$ a prime
  power, in which case the requirement that everything is symmetric
  with respect to the whole symmetric group $\sym(k)$ can be weakened
  to symmetry with respect to the subgroup $G=(\mathbb Z_p)^r$.
\end{cor}
\begin{proof}
  As stated in Section \ref{sec:bu-tuck-prop}, $G$ has the
  Tucker-property. Now apply the previous proposition.
\end{proof}


\section{Towards constructive proofs in topological combinatorics}

The progress in further pursuing the scheme discussed in
Section~\ref{sec:matousek-proof} heavily depends on the question
whether there is a \emph{constructive proof of the Tucker-property
for $G$} at least in the case where $G=\mathbb Z_p$ for $p>2$ a
prime. Such a proof would dramatically increase the chances for
obtaining combinatorial proofs for many theorems from topological
combinatorics.\par

The scheme for discovering constructive proofs for combinatorial
statements, originally deduced by topological arguments, outlined
in Section~\ref{sec:matousek-proof}, is certainly not unique.
There ought to exist other approaches which may be more suitable
for some applications. For example the Tucker lemma and its
generalizations can be incorporated \cite{Sarkaria-90} into a
general problem of finding combinatorial formulas for
Stiefel-Whitney and other characteristic cohomological classes.
More generally, the topological methods used in combinatorial
applications are often naturally seen as part of topological {\em
obstruction theory}. Consequently, there ought to be a close
relationship between the problem of finding constructive proofs
with the program of developing effective obstruction theory which
has been recognized \cite{B-V-Z} as one of the problems paradigmatic
for {\em computational topology} \cite{elisat}.

\subsection*{\bf \textup Acknowledgements}  The authors would like to
thank Francis Su for bringing the consensus-$\frac{1}{k}$-division
problem to their attention, and to G\"unter Ziegler, Arnold
Wassmer and Marc Pfetsch for valuable remarks and suggestions.

\end{document}